\documentclass{amsart}

\usepackage{amscd,amssymb}

\newtheorem{theorem}{Theorem}
\newcommand{\bt}{\begin{theorem}}
\newcommand{\et}{\end{theorem}}

\newtheorem{lemma}{Lemma}
\newcommand{\bl}{\begin{lemma}}
\newcommand{\el}{\end{lemma}}

\newtheorem{corollary}{Corollary}
\newcommand{\bc}{\begin{corollary}}
\newcommand{\ec}{\end{corollary}}

\newcommand{\beq}{\begin{equation}}
\newcommand{\eeq}{\end{equation}}
\newcommand{\benum}{\begin{enumerate}}
\newcommand{\eenum}{\end{enumerate}}
\newcommand{\mcn}{\ensuremath{ \mathcal N}}
\newcommand{\N}{\ensuremath{ \mathbf N }}
\newcommand{\R}{\ensuremath{\mathbf R}}
\newcommand{\T}{\ensuremath{\mathbf T }}
\newcommand{\Z}{\ensuremath{\mathbf Z}}
\DeclareMathOperator{\Add}{Add}
\DeclareMathOperator{\card}{\text{card}}
\DeclareMathOperator{\Geo}{Geo}
\DeclareMathOperator{\id}{id}

\title[Number theory and geometric group theory]{An inverse problem in number theory and geometric group theory}
\author{Melvyn B. Nathanson}
\address{Department of Mathematics\\
Lehman College (CUNY)\\
Bronx, New York 10468}
\email{melvyn.nathanson@lehman.cuny.edu}

\thanks{Supported in part by a grant from the PSC-CUNY Research Award Program.  This paper was written while the author was a visiting fellow in the mathematics department at Princeton University.}

\date{\today}

\begin{document}

\begin{abstract}
This paper describes a new link between combinatorial number theory and geometry.  The main result states that $A$ is a finite set of relatively prime positive integers if and only if $A = (K-K) \cap \N$, where $K$ is a compact set of real numbers such that for every $x\in \R$ there exists $y \in K$ with $x\equiv y \pmod{1}$.  In one direction, given a finite set $A$ of relatively prime positive integers, the proof constructs an appropriate compact set $K$ such that $A = (K-K) \cap \N$.  In the other direction, a strong form of a fundamental result in geometric group theory is applied to prove that $(K-K) \cap \N$ is  a finite set of relatively prime positive integers if $K$ satisfies the appropriate geometrical conditions.   Some related results and open problems are also discussed.
\end{abstract}

\maketitle

\section{From compact sets to integers}
The object of this note is to describe a new connection between number theory and geometry.  Let \R, \Z, and \N\ denote the real numbers, integers, and positive integers, respectively.  For every $x \in \R$, let $[x] \in \Z$ and $(x) \in [0,1)$ denote the integer part and fractional part of $x$.   Let $\Z^n$ denote the additive group of $n$-dimensional lattice points in the Euclidean space $\R^n$.

We recall the following definitions.  
The set $A$ of integers is \emph{relatively prime}, denoted $\gcd(A) = 1$,  if $A$ is nonempty and the elements of $A$ have no common factor greater than 1.  Equivalently, $A$ is relatively prime if $A$ generates the additive group \Z.  
The set $A$ of $n$-dimensional lattice points is relatively prime if the elements of $A$ generate the additive group $\Z^n$.

Let $H$ be a subgroup of a multiplicative group $G$, and let $x$ and $y$ be elements of $G$.  We say that $x$ and $y$ are \emph{congruent modulo $H$}, denoted $x\equiv y \pmod{H}$,  if $xy^{-1} \in H$.  If the group $G$ is additive, then $x\equiv y \pmod{H}$ if $x - y \in H$.  For example, let $G = \R$ and $H = \Z$.  The real numbers $x$ and $y$ are congruent modulo \Z, that is, $x\equiv y \pmod{\Z}$ or, in more traditional notation, $x\equiv y \pmod{1}$, if and only if they have the same fractional part.  

In a multiplicative group $G$ with identity $e$, the \emph{difference set} of a nonempty subset $K$ of $G$ is
\[
KK^{-1} = \{xy^{-1} : x,y \in K \}.
\]
In an additive abelian group $G$, the \emph{difference set} of a subset $K$ of $G$ is
\[
K-K = \{x-y : x,y \in K \}.
\]
Note that $e \in KK^{-1}$ and that the difference set  is symmetric: $z\in KK^{-1}$ if and only if $z^{-1} \in KK^{-1}$ (resp. $z\in K-K$ if and only if $-z \in K-K$).

We are interested in sets of integers contained in difference sets of sets of real numbers.   Our main theorem gives a geometric condition for a finite set of positive integers to be relatively prime.  The geometry uses the concept of an $\mathcal{N}$-set, which is a compact subset $K$ of $\R^n$ such that for every $x\in \R^n$ there exists $y \in K$ with $x\equiv y \pmod{\Z^n}$.

\bt   \label{IntegersCompact:theorem:integers}
Let $A$ be a finite set of positive integers.  The set $A$ is relatively prime if and only if there exists an $\mathcal{N}$-set $K$ in \R\ such that $A = (K-K)\cap \N$.
\et

In Section~\ref{IntegersCompact:section:InverseProblem}, we solve the inverse problem:  Given a finite set of relatively prime integers, we construct an  $\mathcal{N}$-set $K$ in \R\ such that $A = (K-K)\cap \N$.  In Section~\ref{IntegersCompact:section:ApplyGGT} we  prove that if  $K$ is an $\mathcal{N}$-set in $\R$, then $A = (K-K)\cap \N$ is a finite set  of relatively prime  positive integers.   The proof uses the ``fundamental observation of geometric group theory'' (Theorem~\ref{IntegersCompact:theorem:FundamentalObservation}), which is reviewed in Section~\ref{IntegersCompact:Appendix}.

Ideas from geometric group theory have been used recently to obtain new results in number theory (e.g. Nathanson~\cite{nath09b,nath09d,nath09c}), and should continue to be useful.  The book of  de la Harpe, \emph{Topics in Geometric Group Theory}~\cite{harp00}, is an excellent survey of this subject.  Theorem~\ref{IntegersCompact:theorem:FundamentalObservation} was discovered and proved independently by Efremovi{\v{c}}~\cite{efre53}, {\v{S}}varc~\cite{svar55}, and Milnor~\cite{miln68a}.

\section{The inverse problem}      \label{IntegersCompact:section:InverseProblem}
In this section we prove that every finite set of relatively prime positive integers can be realized as the difference set of an $\mathcal{N}$-set.  The construction depends on the following simple observation.

\bl
Let $K$ be a set of real numbers, and let $a \in \Z \setminus \{ 0\}$.  Then $a \in K-K$ if and only if there is a two-element subset $\{x,y\}$ of $K$ such that $(x) =  (y)$ and $a = [x] - [y]$.
\el

\begin{proof}
For any non-zero integer $a$, we have $a \in K-K$ if and only if there exist $x,y \in K$ such that $x\neq y$ and 
\[
a = x-y = [x] -[y] +(x) - (y).
\]
Since $[x]-[y] \in \Z$ and $-1 < (x)-(y) < 1$, it follows that $(x)-(y)=0$ and $a = [x]-[y]$.  The set $\{x,y\}$ satisfies the conditions of the Lemma.  
\end{proof}

Here are three examples.  We associate the set $A_1 = \{ 2, 5\}$ with the $\mathcal{N}$-set 
\[
K(A_1) = [0, 1/3] \cup [2+1/3, 2+2/3] \cup [4+2/3,5].
\]
There are three two-element subsets $\{x,y\}$ of $K(A_1)$ such that $x$ and $y$ have the same fractional part:  $\{1/3,2+1/3\}$, $\{2+2/3,4+2/3\}$, and $\{0,5\}$.

The set $A_2 = \{ 6,10,15\}$ arises from the $\mathcal{N}$-set 
\[
K(A_2) = [0, 1/3] \cup [15+ 1/3,15+2/3] \cup [9 + 2/3,10].
\]
The complete list of the two-element subsets $\{x,y\}$ of $K(A_2)$ such that $x$ and $y$ have the same fractional part is: $\{1/3,15+1/3\}$, $\{15+2/3,9+2/3\}$, and $\{0,10\}$.

For the set $A_3 = \{18,28,63\}$, the $\mathcal{N}$-set 
\begin{align*}
K(A_3) = \bigcup_{i=0}^9 & [-18i+i/13, -18i+(i+1)/13] \cup [-99+10/13,-99+11/13] \\
& \cup  [-36+11/13,-36+12/13] \cup  [27+12/13,28]
\end{align*}
satisfies
\[
\left( K(A_3)  - K(A_3) \right) \cap \N = A_3.
\]
There are exactly 12 two-element subsets $\{x,y\}$ of $K(A_3)$ such that  $x$ and $y$ have the same fractional part.

 In the following Lemma we construct an important example of an $\mathcal{N}$-set on the real line, and its associated difference set of integers.

\bl          \label{IntegersCompact:lemma:FundamentalExample}
For the positive integer $w$, let 
\[
\lambda_0 < \lambda_1 < \cdots < \lambda_{w-1}  < \lambda_w 
\]
be a strictly increasing sequence of real numbers such that 
\[
\lambda_w = \lambda_0+1
\]
and let $b_0, b_1, \ldots,  b_{w-1}$ be a sequence of integers such that 
\[
b_{k-1} \neq b_k \text{ for $k = 1,\ldots, w-1$} 
\]
and
\[
1 + b_{w-1} \neq b_0.
\]
The set
\[
K' = \bigcup_{k=0}^{w-1} [b_k + \lambda_k,b_k+\lambda_{k+1}]
\]
is an $\mathcal{N}$-set, and 
\[
(K' - K') \cap \N =  \{ |b_k - b_{k-1}| : k = 1,\ldots, w-1 \} \cup \{ |1 + b_{w-1} - b_0| \}
\]
is a finite set of relatively prime positive integers.
\el

\begin{proof}
The set $K'$ is compact because it is a finite union of closed intervals, and an $\mathcal{N}$-set because
\[
\bigcup_{k=0}^{w-1}[\lambda_k,\lambda_{k+1}] = [\lambda_0,\lambda_w] = [\lambda_0,\lambda_0 + 1].
\]
Let $A$ be the finite set  of positive integers contained in the difference set $K'-K'$.  Since
\[
\{ \{ b_{k-1}+\lambda_k,b_k+\lambda_k \} : k=1,\ldots, w-1\} \cup 
\{ \{ b_0+ \lambda_0, b_{w-1}+\lambda_w\} \}
\]
is the set of all two-element subsets $\{ x,y \}$ of $K'$ with $(x) = (y)$, it follows that 
\[
A = (K'-K') \cap \N =  \{ |b_k - b_{k-1}| : k = 1,\ldots, w-1 \} \cup \{ | 1 + b_{w-1} - b_0 | \}.
\]
Choose $\varepsilon_k \in \{1,-1\}$ such that 
\[
|b_k - b_{k-1}| = \varepsilon_k  (b_k - b_{k-1})
\]
for $k=1,\ldots, w-1$, and $\varepsilon_w \in \{1,-1\}$ such that 
\[
| 1 + b_{w-1} - b_0| = \varepsilon_w ( 1 + b_{w -1} - b_0).
\]
Since
\[
1 = \varepsilon_w | 1 + b_{w-1}-b_0 | - \sum_{k=1}^{w-1} \varepsilon_k | b_k - b_{k-1}|
\]
it follows that $A$ is a finite set of relatively prime positive integers.  
\end{proof}

\bt        \label{IntegersCompact:theorem:InverseConstruction}
If $A$ is a finite set of relatively prime positive integers, then there is an $\mathcal{N}$-set $K$ such that  $A = (K-K) \cap \N$.
\et

\begin{proof}
Since the elements of $A$ are relatively prime, we can write 1 as an integral linear combination of elements of $A$.   Thus, there exist pairwise distinct integers $a_1,\ldots, a_h$ in $A$, positive integers $w_1,\ldots, w_h$, and $\varepsilon_1,\ldots, \varepsilon_h \in \{ 1, -1\}$  such that 
\beq   \label{IntegersCompact:Represent1}
\sum_{i=1}^h \varepsilon_i w_i a_i =1.
\eeq 
Rewriting~\eqref{IntegersCompact:Represent1}, we obtain 
\beq   \label{IntegersCompact:Represent2}
\varepsilon_h a_h = 1 + \sum_{i=1}^{h-1} w_i(-\varepsilon_i a_i )+ (w_h -1)(-\varepsilon_h a_h).
\eeq
Let $w_0 = 0$ and  $w = \sum_{i=1}^h w_i.$   For $j = 1,2,\ldots, w$, we define integers $\tilde{a}_j$ as follows:
If
\[
w_1+\cdots + w_{i-1}+1 \leq j \leq  w_1+\cdots + w_{i-1}+ w_i
\]
then 
\[
\tilde{a}_j = -\varepsilon_i a_i.
\]
It follows that
\[
1 + \sum_{j=1}^w \tilde{a}_j = 1 + \sum_{i=1}^h w_i(-\varepsilon_i a_i) = 0.
\]
For $k = 0,1,\ldots, w$, we consider the  integers 
\beq   \label{CompactInteger:bk}
b_k =  \sum_{j=1}^{k} \tilde{a}_j
\eeq
and real numbers 
\[
\lambda_k = \frac{k}{w}.
\]
Then  $b_0 = 0$, 
\[
0 = \lambda_0 < \lambda_1 < \cdots < \lambda_w = 1
\]  
and, for $k=1,\ldots, w$,
\[
b_k - b_{k-1} = \tilde{a}_k \neq 0.
\] 
It follows from~\eqref{IntegersCompact:Represent1} and~\eqref{CompactInteger:bk} that 
\begin{align*}
1 + b_{w-1} 
& =  1 + \sum_{j=1}^{w-1} \tilde{a}_j =  1 + \sum_{j=1}^{w} \tilde{a}_j  - \tilde{a}_w\\  
& =  -\tilde{a}_w = \varepsilon_h a_h \neq 0 = b_0.
\end{align*}
Construct the $\mathcal{N}$-set 
\[
K' = \bigcup_{k=0}^{w-1} [b_k+\lambda_k, b_k+\lambda_{k+1}].
\]
By Lemma~\ref{IntegersCompact:lemma:FundamentalExample},
\begin{align*}
(K'-K') \cap \N & =  \{ |b_k - b_{k-1}| : k = 1,\ldots, w-1 \} \cup \{ |1 + b_{w-1}-b_0 | \} \\
& =  \{ |\tilde{a}_k| : k = 1,\ldots, w-1 \} \cup \{ |-\tilde{a}_w | \} \\
& =  \{a_i : i = 1,\ldots, h \}.
\end{align*}
Let $\card(A) = \ell.$  If $\ell = h$, then $A =  \{ a_1,\ldots, a_h \}$ and we  set  $K = K'$.  

Suppose that  $\ell > h$ and $A \setminus \{ a_1,\ldots, a_h \} = \{ a_{h+1},\ldots, a_{\ell} \} \neq \emptyset$.  Since 
\[
\frac{i}{w(\ell - h + 1)} \in \left[ 0 , \frac{1}{w} \right] = [b_0+\lambda_0,b_0+\lambda_1] \subseteq K'
\] 
for $i = 1,2,\ldots, \ell -h$, it follows that 
\[
K = K' \cup \left\{ a_{h+i}+\frac{i}{w(\ell - h + 1)} : i = 1,2,\ldots, \ell -h \right\}
\]
is an $\mathcal{N}$-set such that $A = (K-K)\cap \N$.  This completes the proof.
\end{proof}

Let $A$ be a finite set of relatively prime positive integers.  
We define the weight of a representation of 1 in the 
form~\eqref{IntegersCompact:Represent1} by
\[
\sum_{i=1}^h w_i + \card(A) - h.
\]
We define the \emph{additive weight} of $A$, denoted $\Add(A)$ as the smallest weight of a representation of 1 by elements of $A$.  Note that $\Add(A) \geq \card(A)$ for all $A$, and $\Add(A) = \card(A)$ if and only if there exist distinct integers $a_1,\ldots, a_h \in A$ and $\varepsilon_1, \ldots, \varepsilon_h \in \{ 1, -1\}$ such that $\sum_{i=1}^h \varepsilon_i a_i = 1$.

We define the weight of an \mcn-set $K$ as the number of connected components of $K$, and the 
\emph{geometric weight} of $A$, denoted $\Geo(A)$ as the smallest weight of an \mcn-set $K$ such that $A = (K-K)\cap \N$.  

The following result follows immediately from the proof of 
Theorem~\ref{IntegersCompact:theorem:InverseConstruction}.

\bc
Let $A$ be a finite set of relatively prime positive integers.  Then
\[
\Geo(A) \leq \Add(A).
\]
\ec

There exist sets $A$ such that $\Geo(A) < \Add(A).$  For example, if $A = \{1,2,3,\ldots, n\}$, then $K = [0,n]$ is an \mcn-set of weight 1 such that $(K-K)\cap \N= A$, and so $\Geo(A) = 1 < n = \Add(A)$.

\section{Relatively prime sets of lattice points} 
\label{IntegersCompact:section:ApplyGGT}
In this section we obtain the converse of Theorem~\ref{IntegersCompact:theorem:InverseConstruction}.

\bt        \label{IntegersCompact:theorem:converse}
If $K$ is an $\mathcal{N}$-set in $\R$, then $A = (K-K)\cap \N$ is a finite set  of relatively prime  positive integers.
\et

We prove this result in $n$ dimensions.

\bt        \label{IntegersCompact:theorem:latticepoints}
If $K$ is an $\mathcal{N}$-set in $\R^n$, then $A = (K-K)\cap \Z^n$ is a finite set  of relatively prime lattice points.
\et

Note the necessity of the compactness condition.  For $n \geq 1$, the non-compact set $K = [0,1)^n$ has the property that for all $x\in \R^n$ there exists $y \in K$ with $x \equiv y \pmod{\Z^n}$, but $(K-K)\cap \Z^n = \{0\}$.

\begin{proof}
The proof uses a result called  ``the fundamental observation of geometric group theory'' (Theorem~\ref{IntegersCompact:theorem:FundamentalObservation}). We discuss this in Appendix~\ref{IntegersCompact:Appendix}.

The additive group $\Z^n$ acts isometrically and properly discontinuously on $\R^n$ by translation: $(g,x) \mapsto g+x$ for $g \in \Z^n$ and $x \in \R^n$.  The quotient space $\Z^n \setminus \R^n$ is the $n$-dimensional torus, which is compact, and so the group action $\Z^n \curvearrowright \R^n$ is  co-compact.  
Let $\pi:\R^n \rightarrow \Z^n \setminus \R^n$ be the quotient map.  Then $\pi(x) = \langle x \rangle$ is the orbit of $x$ for all $x \in \R^n$.  
If $K$ is an $\mathcal{N}$-set in $\R^n$, then $K$ is compact, and for every $x\in \R^n$ there exists $y \in K$ such that $x \equiv y \pmod{\Z^n}$.  This means that $\pi(y) = \langle x \rangle$, and so $\pi(K) = \Z^n \setminus \R^n$.  Applying Theorem~\ref{IntegersCompact:theorem:FundamentalObservation} to the set $K$, we conclude that the set 
\[
A = \{ a \in \Z^n : K \cap (a+K) \neq \emptyset\}
\]
is a finite set of generators for $\Z^n$.   Moreover, $a\in A$ if and only if $a\in \Z^n$ and there exists $x \in K$ such that  $x \in a+K$, that is, $x = a+y$ for some $y \in K$.  Equivalently, $a\in A$ if and only if $a = x-y \in (K-K) \cap \Z^n$.  This proves Theorem~\ref{IntegersCompact:theorem:latticepoints}.
\end{proof}

The symmetry of the difference set immediately implies Theorem~\ref{IntegersCompact:theorem:converse}.

We can state the following general inverse problem in geometric group theory:  If $A$ is a finite set of generators for a group $G$ such that $A$ is symmetric and contains the identity of $G$,  does there exist a geometric action of $G$ on a  metric space $X$ with quotient map $\pi: X \rightarrow G\setminus X$ such that 
$A = \{ a \in G : K \cap aK \neq \emptyset\}$ for some compact set $K$ with $\pi(K) = G\setminus X$?  
If $X$ is a group and $G$ is a subgroup of $X$ that acts on $X$ by left translation, then 
\[
\{ a \in G : K \cap aK \neq \emptyset\} = KK^{-1} \cap G.
\]
If $G = \Z^n$ and $X = \R^n$, then the inverse problem is to determine if every finite symmetric set relatively prime lattice points that contains 0 is of the form $(K-K)\cap \Z^n$ for some $\mathcal{N}$-set $K$.  
In this paper we proved that the answer is ``yes'' for $G = \Z$, but the answer is not known for higher dimension.  In particular, the inverse problem for lattice points is open for $n=2$.  We would like a description of the sets of lattice points that can be represented in the form $(K-K)\cap \Z^2$ for some $\mathcal{N}$-set $K$.

\appendix
\section{The fundamental observation of geometric group theory}
 \label{IntegersCompact:Appendix}

The proof of Theorem~\ref{IntegersCompact:theorem:latticepoints} is an application of what is often called the ``fundamental observation of geometric group theory''~\cite[Chapter IV, pp. 87--88]{harp00}.  We shall describe this result, which is not well known to number theorists.  

We begin by introducing the class of boundedly compact geodesic metric spaces.  The Heine-Borel theorem \index{Heine-Borel theorem} states that, in Euclidean space $\R^n$ with the usual metric, a closed and bounded set is compact.  We shall call a metric space $(X,d)$ \emph{boundedly compact } \index{boundedly compact  metric space}\index{metric space!boundedly compact } if every closed and bounded subset of $X$ is compact.  Equivalently, $X$ is boundedly compact if every closed ball 
\[
B^*(x_0,r) = \{x \in X: d(x_0,x) \leq r\}
\]
 is compact for all $x_0 \in X$ and $r \geq 0$.  
Boundedly compact metric spaces are also called \emph{proper} metric spaces.

A metric space $(X,d)$ is \emph{geodesic}\index{geodesic metric space}\index{metric space!geodesic}  if, for all points $x_0,x_1 \in X$ with $x_0 \neq x_1$,  there is an isometry $\gamma$ from an interval $[a,b] \subseteq \R$ into $X$ such that $\gamma(a) = x_0$ and $\gamma(b) = x_1$.  Thus, if $t, t' \in [a,b]$, then $d(\gamma(t),\gamma(t')) = |t - t'|$.  In particular, $d(x_0,x_1) = d(\gamma(a),\gamma(b)) = b-a$.  For example, let $x_0,x_1 \in \R^n$ with $|x_1-x_0| = T$.  Define $\gamma: [0,T] \rightarrow \R^n$ by
\[
\gamma(t) = x_0 + \frac{t}{T}(x_1-x_0).
\]
Then $\gamma(0)=x_0$, $\gamma(T) = x_1$, and 
\begin{align*}
|\gamma(t)-\gamma(t')| 
& = \left| \left( x_0 + \frac{t}{T}(x_1-x_0) \right)  - \left( x_0 + \frac{t'}{T}(x_1-x_0) \right) \right|  \\
& = \left| \left(\frac{t-t'}{T}\right) ( x_1-x_0)  \right|  = |t-t'|.
\end{align*}
Thus, $\R^n$ is a boundedly compact geodesic metric space.

Let $G$ be a group that acts on a metric space $(X,d)$.  
We say that the group $G$  \emph{acts isometrically} on $X$ if the function $x \mapsto gx$ is an isometry on $X$ for every $g \in G$. 
The  group action is called \emph{properly discontinuous} \index{properly discontinuous group action} if, for every compact subset $K$ of $X$, there are only finitely many $a \in G$ such that $K \cap aK \neq \emptyset.$   Let $A = \{ a \in G : K \cap aK \neq \emptyset\}$.  Then $A \neq \emptyset$ because $e\in A$.  Since
\[
K \cap a^{-1}K = a^{-1} ( K \cap aK)
\]
it follows that $A^{-1} = A$.  

For every element $x_0 \in X$, the \emph{orbit}\index{orbit} of $x_0$ is the set
\[
\langle x_0 \rangle = \{gx_0:g\in G\} = Gx_0.
\]
The orbits of elements of $X$ partition the set $X$.  Let $G\setminus X$ denote the set of orbits of the group action, and define the function  $\pi: X \rightarrow G\setminus X$ by $\pi(x) = \langle x \rangle$.  We call $G\setminus X$ the \emph{quotient space}\index{quotient space} of $X$ by $G$, and we call $\pi$ the \index{quotient map} \emph{quotient map} of $X$ onto $G\setminus X$.  Note that every orbit $\langle x \rangle$ is a subset of the set $X$ and a point in the quotient space $G\setminus X$.  

We define the quotient topology on $G\setminus X$ as follows:  A  set $V$  in $G\setminus X$ is open if and only if $\pi^{-1}(V)$ is open in $X$.  This is the largest topology on the quotient space $G\setminus X$ such that the quotient map $\pi$ is continuous.  
We call the group action $G \curvearrowright X$ \emph{co-compact} \index{co-compact} if the quotient space $G\setminus X$ is compact.   
An isometric, properly discontinuous, co-compact action of a group $G$ on a boundedly compact geodesic metric space is called a \emph{geometric action}.

We now state the ``fundamental observation of geometric group theory.''

\bt   \label{IntegersCompact:theorem:FundamentalObservation}
Let $(X,d) $ be a boundedly compact geodesic metric space and let $G$ be a group that acts isometrically  on $X$.  Suppose that the group action $G \curvearrowright X$ is properly discontinuous and co-compact.   Let $\pi: X \rightarrow G\setminus X$ be the quotient map, and let $K$ be a compact subset of $X$ such that $\pi(K) = G\setminus X$.  Then 
\[
A = \{ a \in G : K \cap aK \neq \emptyset\}
\]
is a finite set of generators for $G$.  
\et

For example, the additive group $\Z^n$ of $n$-dimensional lattice points acts on Euclidean space $\R^n$ by translation:   $\alpha_g(x) =  g+x$ for $g \in \Z^n$ and $x \in \Z^n$.
The group $\Z^n$ acts isometrically on $\R^n$ since
\[
 |\alpha_g(x) - \alpha_g(y)|  = |(g+x) - (g+y)| = |x-y|.
\]

Let $K$ be a compact subset of $\R^n$.  Then $K$ is bounded and there is a number $r>0$ such that $|x| < r$ for all $x\in K$.  If $g\in \Z^n$ and $K \cap (g+K) \neq \emptyset$, then there exists $x \in K$ such that $g+x \in K$.  Therefore, 
\[
|g| - r < |g| - |x| \leq |g+x| < r
\]
and $|g| < 2r$.   There exist only finitely many lattice points in $\Z^n$ of length less than $2r$, and so the action on $\Z^n$ on $\R^n$ is properly discontinuous.  

We shall prove that the group action $\Z^n \curvearrowright \R^n$ is co-compact.  Let $\pi: \R^n \rightarrow \Z^n \setminus \R^n$ be the quotient map.  The quotient space $\T^n = \Z^n \setminus \R^n$ is called the\index{torus} \emph{$n$-dimensional torus}.  Let $\{W_i\}_{i\in I}$ be an open cover of $\T^n$, and define $V_i = \pi^{-1}(W_i)$ for all $i \in I$.  Then $\{V_i\}_{i\in I}$ is an open cover of $\R^n$.   The \emph{unit cube}\index{unit cube} 
\[
K = [0,1]^n = \{ x = (x_1,\ldots, x_n) \in \R^n : 0 \leq x_i \leq 1 \text{ for all } i = 1,\ldots, n\}
\]
is a compact subset of $\R^n$, and $\pi(K) = \T^n$.  Since $\{V_i\}_{i\in I}$ is an open cover of $K$, it follows that there is a finite subset $J$ of $I$ such that $K \subseteq \bigcup_{j\in J} V_j$, and so
\[
\T^n = \pi(K) \subseteq  \bigcup_{j \in J} \pi(V_j) = \bigcup_{j \in J}W_j.
\]
Therefore, $\T^n$ is compact and the group action $\Z^n \curvearrowright \R^n$ is co-compact.

\def\cprime{$'$} \def\cprime{$'$} \def\cprime{$'$}
\providecommand{\bysame}{\leavevmode\hbox to3em{\hrulefill}\thinspace}
\providecommand{\MR}{\relax\ifhmode\unskip\space\fi MR }
\providecommand{\MRhref}[2]{%
  \href{http://www.ams.org/mathscinet-getitem?mr=#1}{#2}
}
\providecommand{\href}[2]{#2}

\end{document}